\documentclass{amsart}
\usepackage{amssymb} 
\usepackage{enumitem}
\usepackage{multirow}
\setlist[itemize]{leftmargin=12mm}
\setlist[enumerate]{leftmargin=12mm}
\usepackage{comment}
\usepackage{url}

\DeclareMathOperator{\Aut}{Aut}

\DeclareMathOperator{\GL}{GL}

\DeclareMathOperator{\SL}{SL}

\DeclareMathOperator{\cond}{cond}

\DeclareMathOperator{\Gal}{Gal}
\DeclareMathOperator{\norm}{Norm}
\DeclareMathOperator{\ord}{\upsilon}
\newcommand{\Q}{{\mathbb Q}}

\newcommand{\Cs}{C_\mathrm{s}}
\newcommand{\Cns}{C_\mathrm{ns}}
\newcommand{\Nns}{N_\mathrm{ns}}
\newcommand{\Ns}{N_\mathrm{s}}

\newcommand{\F}{{\mathbb F}}
\newcommand{\PP}{{\mathbb P}}

\newcommand{\cM}{\mathcal{M}}
\newcommand{\cN}{\mathcal{N}}

\newcommand{\OO}{{\mathcal O}}

\newcommand{\ff}{\mathfrak{f}}
\newcommand{\gf}{\mathfrak{g}}

\newcommand{\fl}{\mathfrak{l}}
     
\newcommand{\sS}{\mathfrak{S}}
\newcommand{\bs}{\mathbf{s}}

\begin {document}

\newtheorem{thm}{Theorem}
\newtheorem{lem}{Lemma}[section]
\newtheorem{prop}[lem]{Proposition}

\theoremstyle{definition}

\theoremstyle{remark}

\title[Serre's Uniformity]{On Serre's Uniformity Conjecture
for Semistable\\ Elliptic Curves
over Totally Real Fields
}
\author{Samuele Anni and Samir Siksek}

\address{Mathematics Institute\\
	University of Warwick\\
Coventry\\
	CV4 7AL \\
	United Kingdom}
\email{samuele.anni@gmail.com}
\email{samir.siksek@gmail.com}

\date{\today}
\thanks{The authors are supported by EPSRC Programme Grant 
\lq LMF: L-Functions and Modular Forms\rq\  EP/K034383/1.
The second-named author is also supported by an EPSRC Leadership
Fellowship EP/G007268/1.}

\keywords{Elliptic curves, Serre's uniformity,
modularity, Galois representation, level lowering, Hilbert modular forms}
\subjclass[2010]{Primary 11F80, Secondary 11G05, 11F41}

\begin{abstract}
Let $K$ be a totally real field, 
and let $S$ be a finite set of non-archimedean
places of $K$. 
It follows from the work of Merel, Momose and David that there
is a constant $B_{K,S}$ so that if $E$ is an elliptic curve
defined over $K$, semistable outside $S$, then for all 
$p>B_{K,S}$, the representation $\overline{\rho}_{E,p}$
is irreducible. We combine this with modularity and level lowering
to show the existence of
an effectively computable constant $C_{K,S}$,
and an effectively computable set of elliptic curves over $K$
with CM
$E_1,\dotsc,E_n$ such that the following holds. 
If $E$ is an elliptic curve over $K$ 
semistable outside $S$, and $p>C_{K,S}$ is prime, then
either $\overline{\rho}_{E,p}$ is surjective, or 
%$\overline{\rho}_{E,p} \sim \overline{\rho}_{\ff_i,p}$ for some $i=1,\dots,n$, where $C_{K,S}$ is an effectively computable constant (depending only on $K$ and $S$) and  the $\ff_i$ belong to a finite computable set of CM Hilbert eigenforms over $K$ of parallel weight $2$ with $\Q$\--rational eigenvalues. 
$\overline{\rho}_{E,p} \sim \overline{\rho}_{E_i,p}$ for some $i=1,\dots,n$.
%where $C_{K,S}$ is an effectively computable constant (depending only on $K$
%and $S$) and $E_1,\dots, E_n$ is an effectively computable
%set of elliptic curves over $K$ with complex
%multiplication.

\end{abstract}
\maketitle

\section{Introduction}
Let $K$ be a number field.
We write $G_K=\Gal(\overline{K}/K)$ for the
absolute Galois group of $K$. For an elliptic curve $E/K$,
we write $\overline{\rho}_{E,p}$
for the associated representation of $G_K$ on the $p$-torsion of $E$:
\[
\overline{\rho}_{E,p}\; :\; G_K \rightarrow \Aut(E[p]) \cong \GL_2(\F_p).
\]

We recall the following celebrated theorem of Serre.
\begin{thm}[Serre {\cite[Th\'eor\`eme~2]{Serre}}] \label{thm:Serre} 
Let $K$ be a number field and $E$ an elliptic curve over $K$
without CM. Then there is a constant $C_{E, K}$ such that for all $p>C_{E,K}$
the representation $\overline{\rho}_{E,p}$ is surjective.
\end{thm}

Serre's Uniformity Conjecture (originally formulated
by Serre as a question \cite[$\S$~4.3]{Serre} and \cite{SerreCH}) asserts the
existence of a constant $C_K$, depending
only on $K$, such that if $E$ is an elliptic curve
over $K$ without complex multiplication, and $p>C_K$ is a prime, then
 the representation $\overline{\rho}_{E,p}$ is surjective. 
Mazur \cite{Mazur} proved that $\overline{\rho}_{E,p}$
is irreducible for any prime $p>163$ and elliptic curve $E$ over $\Q$.
Recently, Bilu, Parent and Rebolledo \cite{BPR}
 proved, for 
$p\ge 11$, $p \ne 13$, and $E/\Q$ without complex multiplication, 
that the image of
the representation is also not contained in the normalizer of a split Cartan
subgroup of $\GL_2(\F_p)$.

No analogues of the above-mentioned theorems of Mazur and of Bilu, Parent and Rebolledo
are known for elliptic curves over general number fields. The strongest
known result is Merel's Uniform Boundedness Theorem \cite{Merel}, which asserts 
the following: for $d \ge 1$, there is a constant $B_d$ such that
if $E$ is an elliptic curve over a number field $K$ of degree $d$,
and $p>B_d$ is a prime, then $E(K)[p]=0$. A number of irreducibility results
are however known for semistable elliptic curves over number fields, whose
proofs make essential use
of Merel's Theorem. For example, Kraus \cite[Appendix B]{KrausIrred} shows that
if $K$ is a number field that does not contain
the Hilbert class field of an imaginary quadratic field, then
there is a constant $B_K$ such that for a prime $p>B_K$
and a semistable elliptic curve $E/K$, the representation $\overline{\rho}_{E,p}$
is irreducible.

As noted by Serre \cite[Theorem~4]{Mazur}, Mazur's Theorem cited above
implies the following: if $E/\Q$ is a semistable elliptic curve without complex multiplication, then the representation $\overline{\rho}_{E,p}$ is surjective for any prime $p\geq 11$. To motivate our present work, it is appropriate to give a sketch of the argument.
By Mazur's Theorem, we may suppose that $\overline{\rho}_{E,p}$ is irreducible. 
As $\Q$ has a real embedding, $\overline{\rho}_{E,p}$ is therefore absolutely irreducible (e.g.\ \cite[Lemma 5]{Rubin}).
If $\overline{\rho}_{E,p}$
is not surjective, then its image is contained in the normalizer $\Nns$ of non-split Cartan subgroup $\Cns$ 
or the normalizer $\Ns$ of a split Cartan subgroup $\Cs$. In either case, the representation $\overline{\rho}_{E,p}$
induces a quadratic character $\psi~:~G_{\Q} \rightarrow N_{*}/C_{*} \cong \{\pm 1\}$. This character
is unramified away from the archimedean and additive places. As $E$ is semistable, we see that 
$\psi$ is unramified away from $\infty$, and as the narrow class number of $\Q$ is $1$,
we have $\psi=1$. It follows that the image of $\overline{\rho}_{E,p}$ is contained
in $\Cs$ or $\Cns$. These groups are absolutely reducible, giving a contradiction.
Over a number field $K$, the argument breaks down.
First the narrow class number of $K$ maybe greater than $1$. Moreover, 
let $L$ be the narrow class field of $K$. If the image of
$\overline{\rho}_{E,p}$
is contained in the normalizer of a Cartan subgroup, then $\overline{\rho}_{E,p} (G_L)$
is contained in a Cartan subgroup: $\Cs$ or $\Cns$. If the former, then
we can conclude the argument using (say) Kraus' result, provided $L$ does not
contain the Hilbert class field of an imaginary quadratic field. In the latter
case, we do the same if $L$ has some real embedding. 
It is clear, however, that the
argument does not hold in general.

In this paper, we restrict ourselves to totally real fields $K$. 
This allows us to apply modularity and level lowering theorems
to semistable elliptic curves $E/K$ whose mod $p$ image
is contained in the normalizer of a Cartan subgroup. 
\begin{thm}\label{thm:surjective}
Let $K$ be a totally real field, and let $S$ be a finite set of non-archimedean
places of $K$. There are an effectively computable constant $C_{K,S}$, depending
only on $K$ and $S$, and a finite computable set $E_1,\dots, E_n$ of elliptic curves over $K$ with complex multiplication
such that the following holds: if $E$ is an elliptic curve over $K$ 
semistable outside $S$, and $p>C_{K,S}$ is prime, then either
$\overline{\rho}_{E,p}$ is surjective, or $\overline{\rho}_{E,p} \sim
\overline{\rho}_{E_i,p}$ for some $i=1,\dots,n$.
\end{thm}

% We recall the following standard conjecture, which extends the well-known
% Eichler--Shimura Theorem over $\Q$.
% \begin{conj}\label{conj:eich-shi} Let $K$ be a totally real field.
% Let $\ff$ be a Hilbert eigenform over $K$ of parallel weight $2$ with rational Hecke eigenvalues.
% Then there is an elliptic curve $E/K$
% such that Hecke $\mathrm{L}$-function of $\ff$
% is equal to the Hasse--Weil $\mathrm{L}$-function of $E$.
% \end{conj}
% The conjecture is known to hold (c.f.\ \cite[Introduction]{Zhang}) if $K$ has odd degree, or if 
% the level of $\ff$ is not a square.
% The Hilbert modular forms in the previous statement are CM Hilbert eigenforms over $K$ of parallel weight $2$
% with $\Q$\--rational eigenvalues, therefore, as explained in \cite[section 2.2]{Blasius}, they correspond to CM elliptic curves over $K$. 
% The following is an immediate corollary to Theorem~\ref{thm:surjective}.
% \begin{cor} 
% %Assume Conjecture \ref{conj:eich-shi}. 
% Under the hypotheses of Theorem~\ref{thm:surjective}, there are an explicitely computable constant
% $C_{K,S}$ and a computable finite set $E_1,\dots, E_n$ of elliptic curves over $K$ with complex multiplication
% such that the following holds: if $E$ is an elliptic curve over $K$ 
% semistable outside $S$, and $p>C_{K,S}$ is prime, then either
% $\overline{\rho}_{E,p}$ is surjective, or $\overline{\rho}_{E,p} \sim
% \overline{\rho}_{E_i,p}$ for some $i=1,\dots,n$.
% \end{cor}

We would like to thank Fred Diamond, Nuno Freitas, James Newton and
John Voight
for helpful discussions. We would also thank the referee for the valuable suggestions and comments.

\section{Irreducibility of mod $p$ representations of elliptic curves}
To deal with the Borel
images we shall invoke the following theorem due
to Freitas and Siksek \cite{FSirred}, but is in fact a corollary
of the ideas of  
David \cite{DavidI} 
and Momose \cite{Momose}
building on Merel's Uniform Boundedness Theorem~\cite{Merel}.
\begin{thm}[{\cite[Theorem~1]{FSirred}}]\label{thm:irred}
Let $K$ be a totally real Galois number field of degree $d$,
with ring of integers $\OO_K$ and Galois group $G=\Gal(K/\Q)$.
Let $\sS=\{0,12\}^G$, which we think of as the set of sequences of values $0$, $12$
indexed by $\tau \in G$. %Fix an ordering $\tau_1,\tau_2,\ldots,\tau_d$
For $\bs=(s_\tau) \in \sS$ and $\alpha \in K$, define the \textbf{twisted norm associated
to $\bs$} by
\[
\cN_\bs(\alpha)= \prod_{\tau \in G} \tau(\alpha)^{s_\tau}.
\]
Let $\epsilon_1,\dots,\epsilon_{d-1}$
be a basis for the unit group of $K$ (modulo $\pm 1$),
and define
\[
A_\bs:=\norm  \left( \gcd ( ( \cN_\bs(\epsilon_1)-1) \OO_K,\ldots, (\cN_\bs(\epsilon_{d-1})-1  ) \OO_K) \right).
\]
Let $B$ be the least common multiple of the $A_\bs$ taken over all $\bs \ne (0)_{\tau \in G}$,
$(12)_{\tau \in G}$. Then $B \ne 0$. Moreover,
let $p \nmid B$ be a rational prime, unramified in $K$, such that $p \geq 17$ or $p = 11$. If $E/K$ is an elliptic curve semistable at all $\upsilon \mid p$
and $\overline{\rho}_{E,p}$ is reducible then $p<(1+3^{6dh})$,
where $h$ is the class number of $K$.
\end{thm}

\section{Modularity}
Let $K$ be a totally real number field, 
%and write $G_K:=\Gal(\overline{\Q}/K)$ for its absolute Galois group. 
and let $E$ be an elliptic curve over $K$.
Recall that $E$
is \textbf{modular} if
there exists a Hilbert cuspidal eigenform $\ff$
over $K$ of parallel weight $2$, with rational
Hecke eigenvalues, such that the Hasse--Weil
$\mathrm{L}$-function of $E$ is equal to the Hecke $\mathrm{L}$-function
of $\ff$. It is conjectured that all elliptic curves
over totally real fields are modular, and, recently, modularity has been proved for elliptic curves over real quadratic fields, see \cite{FHS}. 

For what follows, we need a suitable modularity lifting theorem.
The following such theorem is derived in \cite{FHS}
as a relatively straightforward consequence of
a deep theorem of
Breuil and Diamond \cite[Th\'{e}or\`{e}me 3.2.2]{BreuilDiamond}, which
builds on the work of
Kisin \cite{Kisin}, Gee \cite{Gee}, and
Barnet-Lamb, Gee and Geraghty \cite{BGG1}, \cite{BGG2}.
\begin{thm}[{\cite[Theorem 2]{FHS}}]
\label{thm:BD}
Let $E$ be an elliptic curve over a totally real number field $K$,
and let $p\ne 2$ be a rational prime. 
Suppose $\overline{\rho}_{E,p}$ is modular in the following sense:
$\overline{\rho}_{E,p} \sim \overline{\rho}_{\ff,\varpi}$
for some Hilbert cuspidal eigenform over $K$ of parallel weight $2$,
where $\varpi \mid p$. 
Suppose moreover that $\overline{\rho}_{E,p} (G_{K(\zeta_p)})$
is absolutely irreducible.
Then $E$ is modular.
\end{thm}

\begin{prop}\label{prop:modular}
Let $K$ be a totally real field. Let $p \ge 7$
be a prime that is unramified in $K$. Suppose
that $E$ is semistable at some prime $\upsilon$ of $K$ above $p$,
and that moreover $\overline{\rho}_{E,p}$ is irreducible
but not surjective. Then $E$ is modular.
\end{prop}
\begin{proof}
Write $G:=\overline{\rho}_{E,p}(G_K)$. 
As 
$\upsilon \mid p$ is unramified, we have $K \cap \Q(\zeta_p)=\Q$,
and so $\det \overline{\rho}_{E,p}=\chi : G_K \rightarrow \F_p^*$
is surjective, where $\chi$ is the mod $p$ cyclotomic character.
By assumption $\overline{\rho}_{E,p}$ is irreducible
but not surjective, and so 
$G$ does not contain $\SL_2(\F_p)$.
It follows
\cite[{\S}2]{Serre} that
$G$ is contained in the normalizer of a Cartan subgroup,
or its projectivization $\mathbb{P} G:=G/(G \cap \F_p^*)$
is isomorphic to $A_4$, $S_4$ or $A_5$. In particular, $G$
does not contain elements of order $p$.

Write $I_\upsilon \subset G_K$ for the inertia subgroup at $\upsilon$. 
As $E$ is semistable at $\upsilon$ and $\upsilon$ is an unramified prime, we
have (using \cite[{\S}1.11, {\S}1.12]{Serre} and the fact
that $G$ does not contain elements of order $p$):
\begin{equation}\label{eqn:inertia}
\overline{\rho}_{E,p} \vert_{I_\upsilon}
\sim
\begin{pmatrix}
\chi & 0 \\
0 & 1 
\end{pmatrix} \qquad \text{or} \qquad
\overline{\rho}_{E,p} \vert_{I_\upsilon}\otimes_{\F_p}\F_{p^2}
\sim
\begin{pmatrix}
\omega  & 0 \\
0 & \omega^p 
\end{pmatrix}; 
\end{equation}
here $\omega$ is a level $2$ fundamental character $I_\upsilon \rightarrow \F_{p^2}^*$.
More precisely, if $E$ has good ordinary or multiplicative
reduction at $\upsilon$ then we are in the first case of \eqref{eqn:inertia}, 
and if 
$E$ has good supersingular reduction
at $\upsilon$ then we are in the second case.
We observe from \eqref{eqn:inertia} that $\PP G$ contains
an element of order $p-1$ or $p+1$. 
Since $p \ge 7$, we see that $\PP G$ is not isomorphic to $A_4$, $S_4$
and $A_5$. It follows that $G$ is contained in the normalizer $N_{*}$ of a
Cartan subgroup $C_{*}$. The representation $\overline{\rho}_{E,p}$
is irreducible, and as $K$ is totally real, $\overline{\rho}_{E,p}$ must be absolutely irreducible
(e.g.\ \cite[Lemma 5]{Rubin}). Thus the image $G$ is contained in $N_*$ but not in $C_*$.
\begin{comment}
Let $\psi: ~G_{K} \rightarrow N_{*}/C_{*} \cong \{\pm 1\}$ be
the character induced by the quotient map. The kernel of $\psi$ is an open
subgroup of $G_K$ of index $2$, so its fixed field $L$ is quadratic and
unramified outside $\cond(\overline{\rho}_{E,p})\cdot p$ by construction. 

Let $G_L{=}\Gal(\overline{\Q}/L)$. The representation $\overline{\rho}_{E,p}$
restricted to $G_L$ is absolutely 
reducible since its image is contained in $C_*$.
Therefore, we have that 
$\overline{\rho}_{E,p}\vert_{G_L}=\tau\oplus\tau^\prime$, where
$\tau$ and $\tau^\prime$ are characters of $G_L$. These characters satisfy the
following relation:  $\tau^\prime=\tau^\sigma$, 
where $\sigma$ is the non trivial
element of $\Gal(L/K)$, and
$\tau^\sigma(\gamma)=\tau(\sigma\gamma\sigma^{{-}1})$ for all $\gamma\in G_L$.
Hence, $\overline{\rho}_{E,p}$ 
is the representation of $G_K$ induced by $\tau$. Thus
$\overline{\rho}_{E,p}$ is modular by automorphic induction, see
\cite[Theorem~5.1]{SG}, and arising from a Hilbert modular form $\ff$.
Moreover, it is possible to twist $\psi$ so that $\ff$ has parallel weight two
and $\overline{\rho}_{\ff,\varpi}$ has cyclotomic determinant, where $\varpi
\mid p$, see \cite[Lemma~5.1]{Khare}.
\end{comment}

Now, as $\overline{\rho}_{E,p}$ has solvable image, we can
view it as a totally odd irreducible
 Artin representation. By a standard 
argument (c.f.\ \cite[Proof of Lemma 4.2]{DF}),
we have  
$\overline{\rho}_{E,p} \sim \overline{\rho}_{\ff,\varpi}$,
for some Hilbert modular form $\ff$ over $K$, of
parallel weight $2$, and $\varpi \mid p$.

By Theorem~\ref{thm:BD}, in order to show that $E$ is modular it is sufficient to show that 
$\overline{\rho}_{E,p}(G_{K(\zeta_p)})$ 
is absolutely irreducible. 
Suppose otherwise. 
It follows \cite[Lemma 4.2]{FHS} 
that $G^+:=G \cap \GL^+_2(\F_p)$ is absolutely reducible,
where $\GL_2^+(\F_p)$ is the subgroup of $\GL_2(\F_p)$ consisting
of matrices with square determinant. 
Suppose $E$ has good ordinary or multiplicative reduction
at $\upsilon$ and so we are in the first case of \eqref{eqn:inertia}.
Let $g$ be a generator of $\F_p^*$. Then,
with a suitable choice of basis for $E[p]$, the image $G$ contains
all matrices of the form 
$A_r:=(\begin{smallmatrix} g^r & 0 \\ 0 & 1 \end{smallmatrix})$;
these
share the eigenvectors 
$\mathbf{u}=(\begin{smallmatrix} 1 \\0 \end{smallmatrix})$, 
$\mathbf{v}=(\begin{smallmatrix} 0 \\ 1 \end{smallmatrix})$. As $G$
is absolutely irreducible, 
it must contain some matrix $B$
whose eigenvectors $\ne \mathbf{u}$, $\mathbf{v}$.
It follows that $G^+$ contains the pair of matrices
$A_2$,
$B A_s$,
where $s=0$ or $1$ according to whether $\det(B)$ is a square
or non-square. It is easy to check that these do not have 
common eigenvectors, contradicting the absolute reducibility of $G^+$.
If $E$ has good supersingular reduction at $\upsilon$ then
we are in the second case of \eqref{eqn:inertia}. 
%It follows 
%that $G$ contains an entire non-split Cartan subgroup. 
%As $K$ is real, $G$ is absolutely irreducible. 
%and so $G$
%must be the normalizer of a non-split Cartan subgroup.
It is now easy to check, similarly to the above, that
$G^+$ is absolutely irreducible, giving a contradiction.
This completes the proof.
\end{proof}

\section{Level lowering}

In this section, $K$ is a totally real field, and $S$ a finite
set of non-archimedean primes of $K$. Moreover, $p\ge 7$
is a rational prime that is unramified in $K$ such that
$\upsilon \notin S$ for all $\upsilon \mid p$.
\begin{lem}\label{lem:unramified}
Let $E$ be an elliptic curve defined over $K$
that is semistable outside $S$. Suppose that
$\overline{\rho}_{E,p}$ is irreducible but not surjective.
Then
\begin{enumerate}
\item[(i)] $\overline{\rho}_{E,p}$ is unramified at all $\mu \notin S$,
$\mu \nmid p$;
\item[(ii)] $\overline{\rho}_{E,p}$ is finite at all $\upsilon \mid p$.
\end{enumerate}
\end{lem}
\begin{proof}
Let $\upsilon \mid p$. We would like to prove (ii), which is
certainly true if $E$ has good reduction at $\upsilon$. By hypothesis, 
$E$ is semistable at $\upsilon$, and so we may assume that $E$
has multiplicative reduction at $\upsilon$. 
Write $G_\upsilon \subset G_K$ for the decomposition group at $\upsilon$.
By the proof of Proposition~\ref{prop:modular}, we know that 
$G=\overline{\rho}_{E,p}(G_K)$ does not contain any elements of 
order $p$. It immediately follows that $\overline{\rho}_{E,p} \vert_{G_\upsilon}$
is \lq\lq peu ramifi\'{e}\rq\rq, proving (ii).

Let $\mu$ be a non-archimedean prime of $K$, not in $S$,
and not above $p$. Then $E$ is semistable at $\mu$,
and so the inertia subgroup $I_\mu \subset G_K$ acts
unipotently on $E[p]$. 
As $G$ does not contain elements of order $p$, we have $\overline{\rho}_{E,p}(I_\mu)=1$,
proving (i).
\end{proof}

Now let
\[
\cM=\prod_{\fl \in S} \fl^{2 + 6\ord_\fl(2)+3 \ord_\fl(3)}.
\]
\begin{lem}\label{lem:ll}
Assume the hypotheses of Lemma~\ref{lem:unramified}.
Then there exists a Hilbert eigenform $\ff$ over $K$
of parallel weight $2$ and level dividing $\cM$
such that $\overline{\rho}_{E,p} \sim \overline{\rho}_{\ff,\varpi}$
where $\varpi \mid p$ is a prime of $\Q_\ff$, the field generated
by the eigenvalues of $\ff$. 
\end{lem}
\begin{proof}
Let $\cN$ be the conductor of $E$. The additive part of $\cN$
divides $\cM$ 
(e.g.\ \cite[Theorem~IV.10.4]{SilvermanII}). By Proposition~\ref{prop:modular} and Theorem~\ref{thm:BD}, 
there is a Hilbert eigenform $\ff_0$ over $K$, with rational eigenvalues, 
level $\cN$ and parallel weight $2$ such that
$\overline{\rho}_{E,p} \sim \overline{\rho}_{\ff_0,p}$.
By Lemma~\ref{lem:unramified}, we have $\overline{\rho}_{\ff_0,p}$
is finite at all $\upsilon \mid p$, and unramified at all
$\mu \nmid \cM$. Applying level lowering 
theorems due 
Fujiwara \cite{Fuj},
Jarvis \cite{Jarv} and 
Rajaei \cite{Raj}, 
we may remove these primes from the level (without changing the weight);
the argument is practically identical to that in \cite[Theorem 7]{FS},
and so we omit the details.
\end{proof}

\noindent \textbf{Remark.}
Chen \cite{Chen} observes that if $E$ is an elliptic curve
over $\Q$, and $\overline{\rho}_{E,p}$ has image 
contained in the
normalizer of a Cartan subgroup, then $p$ is a congruence
prime for the newform attached to $E$. Our Lemma~\ref{lem:ll}
encompasses Chen's observation.

\section{Proof of Theorem~\ref{thm:surjective}}
Assume the hypotheses of Theorem~\ref{thm:surjective}:
in particular, let $E$ be an elliptic curve semistable outside $S$. With the help of Theorem~\ref{thm:irred},
we know that there is an effectively computable constant $C_{K,S}$
such that if $p>C_{K,S}$ then $p$ is unramified in $K$, all the primes
$\upsilon \mid p$ satisfy $\upsilon \notin S$, and $\overline{\rho}_{E,p}$
is irreducible. Suppose $\overline{\rho}_{E,p}$ is not
surjective. We now apply Lemma~\ref{lem:ll} to deduce that 
$\overline{\rho}_{E,p} \sim \overline{\rho}_{\ff,\varpi}$
for some cuspidal Hilbert eigenform of parallel weight $2$ and
level dividing $\cM$. 
There are certainly only finitely many such eigenforms.
We would like to increase $C_{K,S}$
by an effectively computable amount so that the conclusion
of Theorem~\ref{thm:surjective} holds. Crucial to the effectivity
is the existence of an algorithm \cite{DV} for determining
the eigenforms $\ff$ of a given weight and level, as well
as their Hecke eigenvalues at given primes, and the fields
generated by these eigenvalues. 
We will eliminate
all such eigenforms $\ff$ with $\Q_\ff \ne \Q$,
where $\Q_\ff$ is the field generated by the eigenvalues of $\ff$.
So suppose that $\overline{\rho}_{E,p} \sim \overline{\rho}_{\ff,\varpi}$
where $\Q_\ff \ne \Q$. Let $\fl$ be the prime
ideal of smallest possible norm such that $\fl \notin S$
and $a_\fl(\ff) \notin \Q$. If $\fl \mid p$, then
$p \mid \norm_{K/\Q}(\fl)$ and so we obtain a contradiction
by supposing that $C_{K,S} > \norm_{K/\Q}(\fl)$.
We may therefore suppose that $\fl \nmid p$.
Comparing the traces of the images of Frobenius at $\fl$
in the representations $\overline{\rho}_{E,p}$
and $\overline{\rho}_{\ff,\varpi}$ we have either 
$a_\fl(\ff) \equiv a_\fl(E) \pmod{\varpi}$
if $E$ has good reduction at $\fl$, or 
$a_\fl(\ff) \equiv \pm (\norm_{K/\Q}(\fl)+1) \pmod{\varpi}$
if $E$ has multiplicative reduction at $\fl$. In the former
case,
by the Hasse--Weil bounds,
$p$ divides  
\[
\prod_{\lvert t \rvert \le B} \norm_{\Q_\ff/\Q} (a_\fl(\ff)-t),
\qquad 
B=2 (\norm_{K/\Q}(\fl))^{1/2}. 
\]
As $a_\fl(\ff) \notin \Q$, all the terms in the product
are non-zero, and so this gives a bound on $p$.
By taking $C_{K,S}$ larger than this product
we obtain a contradiction. If $E$ has multiplicative
reduction at $\fl$, then $p$ divides
\[
\norm_{\Q_\ff/\Q} (a_\fl(\ff) - \norm_{K/\Q}(\fl)-1)
\cdot
\norm_{\Q_\ff/\Q} (a_\fl(\ff) + \norm_{K/\Q}(\fl)+1)
\]
and again we obtain a contradiction by taking $C_{K,S}$
larger than this product.
Thus we are reduced to finitely many forms $\ff$ satisfying
$\Q_\ff=\Q$. 

So far we proved that there are an effectively computable constant $C_{K,S}$ and a finite computable set
$\ff_1,\dots,\ff_n$ of Hilbert eigenforms over $K$ of parallel weight $2$ with $\Q$\--rational eigenvalues such that the following holds: if $E$ is an elliptic
curve over $K$ semistable outside $S$, and $p>C_{K,S}$ is prime, then either $\overline{\rho}_{E,p}$ is surjective, or $\overline{\rho}_{E,p} \sim \overline{\rho}_{\ff_i,p}$ for some $i=1,\dots,n$. 

%To complete the proof of Theorem~\ref{thm:surjective}, 
Next we have to show that the surviving forms $\ff$ have CM, possibly after
enlarging $C_{K,S}$ by an effective amount. 
In fact,
by a theorem of Dimitrov (\cite[Theorem~2.1]{DimitrovCRM},
\cite[$\S$~3]{Dimitrov}),
if $\ff$ does not have CM, there is a constant $B_\ff$
such that for $p>B_\ff$ and $\varpi \mid p$, 
the image of $\overline{\rho}_{\ff,\varpi}$ contains
a conjugate of $\SL_2(\F_p)$. It is however unclear to us
as to whether Dimitrov's proof can be made effective, and so
we proceed in a more elementary manner. 

By the proof
of Proposition~\ref{prop:modular} the image of $\overline{\rho}_{E,p}$
is dihedral, and so there is a quadratic
character $\psi$
such that $\overline{\rho}_{E,p} \sim \overline{\rho}_{E,p} \otimes \psi$. 
It is immediate from Lemma~\ref{lem:unramified} that $\psi$
is unramified away from $S$, the archimedean primes, and the primes
$\upsilon \mid p$.
%Since the two representations are isomorphic, they have the same conductor. Hence,  
%$\psi$ can ramify only at primes dividing $p$ and the conductor of
%$\rho_{E,p}$. 
Suppose $\upsilon \mid p$.
Comparing the restriction of the representation to the inertia subgroup at $\upsilon$,
displayed in \eqref{eqn:inertia}, and the restriction of the twisted
representation by $\psi$, it is easy to deduce that the quadratic character
$\psi$ is unramified at $\upsilon$.  
%Analogously,
%using Lemma~\ref{lem:unramified}, it follows that $\psi$ is unramified outside
%$S$ and the archimedean primes. 
Hence, its conductor divides
$\prod_{\fl \in S} \fl^{1 + 2\ord_\fl(2)}$ and so $\psi$
belongs to a finite effectively computable set of characters.

Suppose $\overline{\rho}_{E,p} \sim \overline{\rho}_{\ff,p}$ where now $\ff$
has rational eigenvalues. Let $\gf=\ff \otimes\psi$.
If $\gf=\ff$ then $\ff$ has CM as desired. Thus
we may suppose $\gf \ne \ff \otimes \psi$. Let $\fl$
be the prime ideal of $K$ of smallest possible norm
so that $\fl \notin S$ and $a_\fl(\ff) \ne a_\fl(\gf)$.
As before, if $\fl \mid p$ or if $E$ has multiplicative
reduction at $\fl$, then we obtain a bound on $p$.
We therefore suppose that $\fl \nmid p$ 
and $E$ has good reduction at $\fl$.
From the relations $\overline{\rho}_{E,p} \sim \overline{\rho}_{E,p} \otimes \psi$
and $\gf=\ff \otimes \psi$ we have $a_\fl(\ff) \equiv a_\fl(\gf) \pmod{p}$.
As $a_\fl(\ff) \ne a_\fl(\gf)$ we obtain a bound on $p$.

%This completes the proof.

Now as the surviving forms $\ff_i$ are CM
Hilbert eigenforms over $K$ of parallel weight $2$ with $\Q$\--rational
eigenvalues. As explained in \cite[$\S$~2.2]{Blasius}, they correspond to
CM elliptic curves $E_i$ over $K$. The conductors of $E_i$
are the levels of $\ff_i$. As there is an effective algorithm
to determine elliptic curves of a given conductor (c.f.\ \cite{CreLing}), the proof
is complete.

\end{document}